\documentclass{birkart}
\usepackage{amsmath}
\usepackage{amssymb}
\usepackage{amsfonts}
\usepackage{graphicx}
\usepackage{texdraw}
\usepackage{graphpap}
\usepackage{wasysym}
\usepackage{pxfonts}
\usepackage[OT2,T1]{fontenc}
\DeclareSymbolFont{cyrletters}{OT2}{wncyr}{m}{n}
\DeclareMathSymbol{\berd}{\beta}{cyrletters}{"42}

\input txdtools

 \newtheorem{thm}{Theorem}[section]
 \newtheorem{cor}[thm]{Corollary}
 \newtheorem{lem}[thm]{Lemma}
 \newtheorem{prop}[thm]{Proposition}
 \theoremstyle{definition}

 \theoremstyle{remark}
 
 \newtheorem{rem}[thm]{Remark}

\numberwithin{equation}{section}
\numberwithin{figure}{section}

\newcommand{\e}{\mathrm e}

\newcommand{\C}{{\mathbb C}}
\newcommand{\D}{{\mathbb D}}

\newcommand{\R}{{\mathbb R}}

\newcommand{\Z}{{\mathbb Z}}

\newcommand{\id}{\operatorname{\text{\bf I}}}

\newcommand{\group}{{\mathfrak G}(\beta)}

\newcommand{\re}{\operatorname{Re}}
\newcommand{\im}{\operatorname{Im}}

\newcommand{\Tope}{{\mathbf T}}
\newcommand{\Sop}{{\mathbf S}}

\newcommand{\Cop}{{\mathbf C}}

\newcommand{\diff}{{\mathrm d}}
\newcommand{\imag}{{\mathrm i}}

\newcommand{\supp}{\operatorname{supp}}

\newcommand{\eps}{\varepsilon}

\begin{document}
%
\title[Heisenberg uniqueness pairs and the Klein-Gordon equation]
{Heisenberg uniqueness pairs and the Klein-Gordon equation}

\author[Hedenmalm]
{H\aa{}kan Hedenmalm}

\address{Hedenmalm: Department of Mathematics\\
The Royal Institute of Technology\\
S -- 100 44 Stockholm\\
SWEDEN}

\email{haakanh@math.kth.se}

\thanks{Research partially supported by the G\"oran Gustafsson Foundation and
by the Swedish Science Council (Vetenskapsr\aa{}det).}


\author[Montes]
{Alfonso Montes-Rodr\'\i{}guez}

\address{Montes-Rodr\'\i{}guez: Department of Mathematical Analysis\\
University of Sevilla\\
Sevilla\\
SPAIN}

\email{amontes@us.es}

\thanks{Research partially supported by Plan Nacional I+D+I grant no. 
MTM2006-09060 and by Junta de Andaluc\'\i{}a grants nos. FQM-260 and 
FQM06-02225.}




\begin{abstract} 
A Heisenberg uniqueness pair (HUP) is a pair $(\Gamma,\Lambda)$,
where $\Gamma$ is a curve in the plane and $\Lambda$ is a set in the plane, 
with the following property: any bounded Borel measure $\mu$ in the plane 
supported on $\Gamma$, which is absolutely continuous with respect to arc 
length, and whose Fourier transform $\widehat\mu$ vanishes on $\Lambda$, 
must automatically be the zero measure.
We prove that when $\Gamma$ is the hyperbola $x_1x_2=1$, and $\Lambda$ is the
lattice-cross
$$\Lambda=(\alpha\Z\times\{0\})\cup(\{0\}\times\beta\Z),$$
where $\alpha,\beta$ are positive reals, then $(\Gamma,\Lambda)$ is an 
HUP if and only if $\alpha\beta\le1$; in this situation, the Fourier transform
$\widehat\mu$ of the measure solves the one-dimensional Klein-Gordon equation. 
Phrased differently, we show that
$$\e^{\pi\imag \alpha n t},\,\,\e^{\pi\imag\beta n/t},\qquad n\in\Z,$$
span a weak-star dense subspace in $L^\infty(\R)$ if and only if 
$\alpha\beta\le1$.
In order to prove this theorem, some elements of linear fractional theory
and ergodic theory are needed, such as the Birkhoff Ergodic Theorem.
An idea parallel to the one exploited by Makarov and Poltoratski (in the 
context of model subspaces) is also needed. As a consequence, we solve a
problem on the density of algebras generated by two inner functions raised
by Matheson and Stessin.
\end{abstract}

\maketitle

\addtolength{\textheight}{2.2cm}







\section{Introduction}

\noindent\bf Heisenberg uniqueness pairs. \rm
Let $\mu$ be a finite complex-valued Borel measure in the plane $\R^2$,
and associate to it the Fourier transform 
$$\widehat\mu(\xi)=\int_{\R^2}\e^{\pi\imag\langle x,\xi\rangle}
\diff\mu(x),$$
where $x=(x_1,x_2)$ and $\xi=(\xi_1,\xi_2)$, with inner product
$$\langle x,\xi\rangle=x_1\xi_1+x_2\xi_2.$$
The Heisenberg uncertainty principle states that both $\mu$ and $\widehat\mu$
cannot both be too concentrated to a point (see \cite{H} for the original paper
of Heisenberg, and \cite{HJ} for a more general treatment); in particular, 
they cannot both have compact support. 
Here, we shall study a variation on that theme. Let $\Gamma$ be a smooth 
curve in $\R^2$, or, more generally, a
finite disjoint union of smooth curves. Suppose that $\supp\mu\subset\Gamma$,
and that $\mu$ is absolutely continuous with respect to arc length measure on 
$\Gamma$.
Which sets $\Lambda\subset\R^2$ have the property that
$$\widehat\mu|_{\Lambda}=0 \quad\implies\quad\mu=0?$$
If this is the case, we say that $(\Gamma,\Lambda)$ is a {\em Heisenberg
uniqueness pair}. A dual formulation is that $(\Gamma,\Lambda)$ is a 
Heisenberg uniqueness pair if and only if the functions
$$e_\xi(x)=\e^{\pi\imag\langle x,\xi\rangle},\qquad \xi\in\Lambda,$$
span a weak-star dense subspace in $L^\infty(\Gamma)$. This concept of
Heisenberg uniqueness pairs has many features in common with the notion 
of (weakly) mutually annihilating pairs of Borel measurable sets having 
positive area measure, which appears, for instance, in the book by Havin and 
J\"oricke \cite{HJ}.

The properties of the Fourier transform with respect to
translation and multiplication by complex exponentials show that for all
points $x^*,\xi^*\in\R^2$, we have
$$(\Gamma+\{x^*\},\Lambda+\{\xi^*\})\quad \text{is an HUP}\quad
\Longleftrightarrow \quad (\Gamma,\Lambda)\quad \text{is an HUP},
\leqno{\text{(inv-1)}}$$ 
where HUP is short for ``Heisenberg uniqueness pair''.
Likewise, it is also straightforward to see that if $T:\R^2\to\R^2$ is an
invertible linear transformation with adjoint $T^*$, then
$$(T^{-1}(\Gamma),T^*(\Lambda))\quad \text{is an HUP}\quad
\Longleftrightarrow \quad (\Gamma,\Lambda)\quad \text{is an HUP}.
\leqno{\text{(inv-2)}}$$ 
\medskip

\noindent\bf Algebraic curves and partial differential equations. \rm
Algebraic curves $\Gamma$ are of particular interest, because of their connection
to partial differential equations. That connection follows from the observation
that for polynomials $p$ of two variables,
$$p\bigg(\frac{\partial_1}{\pi\imag},\frac{\partial_2}{\pi\imag}\bigg)
\widehat\mu(\xi)=\int_{\R^2}
\e^{\pi\imag\langle x,\xi\rangle}\,p(x_1,x_2)\,
\diff\mu(x),$$
so that if $p$ is real-valued and $\Gamma$ is the locus of the equation
$$p(x_1,x_2)=0,$$
then 
$$p(x_1,x_2)\,\diff\mu(x_1,x_2)=0$$
identically. Therefore $\widehat\mu$ solves the PDE
\begin{equation}
p\bigg(\frac{\partial_1}{\pi\imag},\frac{\partial_2}{\pi\imag}\bigg)
\widehat\mu(\xi)=0
\label{PDE}
\end{equation}
in the plane. In fact, the equation \eqref{PDE} encodes the requirement 
that $\supp\mu\subset\Gamma$. 

\medskip

\noindent\bf Conic sections. \rm
We shall consider the case when $\Gamma$ is a {\em conic section}, that is, 
the locus of a quadratic equation
$$ax_1^2+bx_2^2+cx_1x_2+dx_1+ex_2+f=0,$$
where $a,b,c,d,e,f$ are real constants. As we only consider the case when 
$\Gamma$ is a curve, this leaves us with the following cases: a 
straight line, two parallel straight lines, a cross, an ellipse, a 
parabola, or a hyperbola. 
\medskip

\noindent\bf The line. \rm Let us look
at the line first, as a model example. By the invariance properties 
\text{(inv-1)} and \text{(inv-2)}, we may assume that $\Gamma=\R\times\{0\}$,
the $x_1$-axis. In this case, $\widehat\mu(\xi)$ depends only on $\xi_1$,
and it is easy to see that $(\Gamma,\Lambda)$ is a Heisenberg uniqueness pair 
if and only
if $\pi_1(\Lambda)$, the  orthogonal projection of $\Lambda$ to the 
$\xi_1$-axis, is dense.
\medskip

\noindent\bf Two parallel lines. \rm If $\Gamma$ is the union
of two parallel lines, we may without loss of generality assume that
$$\Gamma=\R\times\{0,1\}.$$
In this case, we see from the example of the line that it is necessary for
$(\Gamma,\Lambda)$ to be a Heisenberg uniqueness pair that $\pi_1(\Lambda)$
be dense. But something more is needed. An absolutely continuous measure $\mu$
on $\Gamma$ may be written in the form
$$\diff\mu(x)=f(x_1)\diff x_1\diff\delta_0(x_2)+
g(x_1)\diff x_1\diff\delta_1(x_2),$$
where $f,g\in L^1(\R)$ ($\delta_y$ denotes the unit point mass at the 
point $y$), so that
$$\widehat\mu(\xi)=\widehat f(\xi_1)+\e^{\pi\imag\xi_2}\widehat g(\xi_1).$$
Next, we split 
$$\pi_1(\Lambda)=\pi_1^a(\Lambda)\cup\pi_1^b(\Lambda),$$ 
where the two sets are disjoint: $t\in\pi_1^a(\Lambda)$ if there are 
two lifted points $\xi=(\xi_1,\xi_2)$ and $\eta=(\eta_1,\eta_2)$ in $\Lambda$, 
with $\xi_1=\eta_1=t$ and $\xi_2-\eta_2\notin2\Z$, whereas 
$t\in\pi_1^b(\Lambda)$ if the latter does not happen. We quickly find that
\begin{equation}
\widehat f(t)=\widehat g(t)=0,\qquad t\in\pi_1^a(\Lambda).
\label{eq-1999}
\end{equation}
On the other hand, for $t\in\pi_1^b(\Lambda)$, the expression 
$\e^{\pi\imag\xi_2}$ is a well-defined function of $\xi_1=t$, where 
$\xi_2$ stands for any of the points with $(\xi_1,\xi_2)\in\Lambda$; we write
$\chi(t)$ for this unimodular function. 
If $E$ is a closed subset of $\R$ and $t_0\in E$, we say that a function
$\varphi:E\to\C$ is {\em locally the Fourier transform of an $L^1(\R)$ function
around $t_0$} provided  that there exists a small open interval 
$I$ around $t_0$ and a function $\psi$ which is the Fourier transform of an 
$L^1(\R)$ function, such that $\psi=\varphi$ on $E\cap I$.
Let $\pi_1^c(\Lambda)$ consist of those points $t_0\in\pi_1^b(\Lambda)$ 
where $\chi:\pi_1^b(\Lambda)\to\C$ is locally the Fourier transform of an 
$L^1(\R)$ function around $t_0$.

\begin{thm}
$(\Gamma,\Lambda)$ is a Heisenberg uniqueness pair if and only if
$\pi_1^a(\Lambda)\cup(\pi_1^b(\Lambda)\setminus\pi_1^c(\Lambda))$ is dense in 
$\R$.
\end{thm}

\begin{proof}
We observe that
\begin{equation}
\widehat f(t)=-\chi(t)\widehat g(t),\qquad t\in\pi_1^b(\Lambda).
\label{eq-2000}
\end{equation}
If $t\in\pi_1^b(\Lambda)\setminus\pi_1^c(\Lambda)$, this is only
possible if $\widehat g(t)=0$, so that
\begin{equation}
\widehat f(t)=\widehat g(t)=0,\qquad 
t\in\pi_1^b(\Lambda)\setminus\pi_1^c(\Lambda).
\label{eq-2001}
\end{equation}
A combination of \eqref{eq-1999} and \eqref{eq-2001} shows that
$f=g=0$ (so that $\mu=0$) if the set
$\pi_1^a(\Lambda)\cup(\pi_1^b(\Lambda)\setminus\pi_1^c(\Lambda))$ is dense in 
$\R$.

As for the other direction, suppose that $\pi_1(\Lambda)$ is dense in $\R$,
while $\pi_1^a(\Lambda)\cup(\pi_1^b(\Lambda)\setminus\pi_1^c(\Lambda))$ fails 
to be dense in $\R$. We then pick a point $t_0\in\R$ such that an open 
interval $J$ around it has empty intersection with 
$$\pi_1^a(\Lambda)\cup(\pi_1^b(\Lambda)\setminus\pi_1^c(\Lambda)).$$
But then $\pi_1^c(\Lambda)\cap J$ is dense in $J$, and $\chi$ is locally the 
Fourier transform of an $L^1(\R)$ function around $t_0$. We thus find a 
function $\chi_1$ which coincides with $\chi$ on some open interval 
$I\subset J$ with $t_0\in I$, while $\chi_1$ is the Fourier 
transform of an $L^1(\R)$ function. Next, we pick $g\in L^1(\R)$ with 
$\widehat g(t_0)\neq0$, such that
$\supp\widehat g\Subset I$, and define $f\in L^1(\R)$ via
$\widehat f=-\chi_1\widehat g$, so that \eqref{eq-2000} holds. 
This gives us a nontrivial measure $\mu$ with the required 
properties, and so $(\Gamma,\Lambda)$ cannot be a Heisenberg uniqueness pair.
\end{proof}
\medskip

\noindent\bf The cross. \rm If $\Gamma$ is a cross, the 
PDE \eqref{PDE} expresses the wave equation.
By the invariance properties \text{(inv-1)} and \text{(inv-2)}, we may restrict
our attention to the case when
$$\Gamma=(\R\times\{0\})\cup(\{0\}\times\R)$$ 
is the union of the two axes. Here, it appears that the characterization of
uniqueness pairs $(\Gamma,\Lambda)$ may get quite complicated. Obviously, 
it is a necessary condition that $\pi_1(\Lambda)$ and $\pi_2(\Lambda)$ be
dense ($\pi_2(\Lambda)$ is the orthogonal projection to the $\xi_2$-axis).
This is far from sufficient, because if $\Lambda$ is contained in a smooth
graph, we may run into trouble. 
For instance, if $\Lambda$ is contained in the diagonal 
$\xi_1=\xi_2$, then we may choose 
$$\diff\mu(x_1,x_2)=f(x_1)\,\diff x_1\diff\delta_0(x_2)-
f(x_2)\,\diff x_2\diff\delta_0(x_1),$$
where $f\in L^1(\R)$, which is supported on $\Gamma$ and nontrivial 
generically, while 
$\widehat\mu(\xi_1,\xi_2)=0$ for $\xi_1=\xi_2$. 
\medskip

\noindent\bf The ellipse. \rm
If $\Gamma$ is an ellipse, the invariance of \text{(inv-1)} and \text{(inv-2)}
allows us to focus on the circle
$$\Gamma=\{x=(x_1,x_2)\in\R^2:\,\,x_1^2+x_2^2=1\}.$$
The corresponding PDE \eqref{PDE} is the eigenvalue equation for the Laplacian.
Here, the fact that $\Gamma$ is compact entails that $\widehat\mu(\xi)$ extends
to an entire function of exponential type in $\C^2$. It would seem that
reasonable criteria on $\Lambda$ may be found that are at least close to being
necessary and sufficient for $(\Gamma,\Lambda)$ to be a Heisenberg uniqueness 
pair. 
\medskip

\noindent\bf The parabola. \rm
If $\Gamma$ is a parabola, the invariance of \text{(inv-1)} and \text{(inv-2)} 
allows us to focus on the parabola
$$\Gamma=\{x=(x_1,x_2)\in\R^2:\,\,x_2=x_1^2\}.$$
The corresponding PDE \eqref{PDE} is the one-dimensional Schr\"odinger 
equation without potential. Here, the problem of characterizing the Heisenberg 
uniqueness pairs $(\Gamma,\Lambda)$ appears quite challenging. 
\medskip

\noindent\bf The hyperbola. \rm
We shall focus most of our attention to the case when $\Gamma$ is a
hyperbola. The corresponding PDE \eqref{PDE} is the one-dimensional 
Klein-Gordon equation.
We will see that the situation with Heisenberg uniqueness pairs is 
dramatically different from that of the cross. By the invariance 
\text{(inv-1)} and \text{(inv-2)}, we may assume that the hyperbola is given 
by
$$x_1x_2=1.$$

\begin{thm}
Suppose $\Gamma$ is the hyperbola $x_1x_2=1$ and that $\Lambda$ is the 
lattice-cross
$$\Lambda=(\alpha\Z\times\{0\})\cup(\{0\}\times\beta\Z),$$
where $\alpha,\beta$ are positive reals. Then $(\Gamma,\Lambda)$ is
a Heisenberg uniqueness pair if and only if $\alpha\beta\le1$.  
\label{thm-1}
\end{thm}

The remainder of this work is devoted to proving this assertion.
But before we turn to the proof, let us consider a generalization
which is more or less immediate.

\begin{cor}
Suppose $\Gamma_\varepsilon$ is the hyperbola $x_1x_2=\varepsilon$, where 
$\varepsilon\neq0$ is real, and that $\Lambda$ is the 
lattice-cross
$$\Lambda=(\alpha\Z\times\{0\})\cup(\{0\}\times\beta\Z),$$
where $\alpha,\beta$ are positive reals. Then $(\Gamma_\varepsilon,\Lambda)$ is
a Heisenberg uniqueness pair if and only if 
$\alpha\beta\le1/|\varepsilon|$.  
\label{cor-1}
\end{cor}

The eccentricity of the hyperbola $\Gamma_\varepsilon$
is $\sqrt{2}$ independently of $\varepsilon$. 
The condition of the corollary ($\alpha\beta\le1/|\varepsilon|$)
gets weaker as $|\varepsilon|$ decreases. However, in the 
limit situation $\varepsilon=0$ -- the cross --
the situation changes dramatically: if $\Lambda$ is contained in the dual 
cross $(\R\times\{0\})\cup(\{0\}\times\R)$, then $\Lambda$ must actually be
dense in the cross for $(\Gamma_0,\Lambda)$ to be a Heisenberg uniqueness pair.

\begin{rem}
\label{rem-1}
Consider for a moment the sets
$$\Lambda'=([\theta,+\infty[\times]-\infty,0])\cup
(]-\infty,0]\times[0,+\infty[)$$
and
$$\Lambda''=\big\{(\xi_1,\xi_2)\in\R^2:\,\,a_1\xi_1+a_2\xi_2=0\big\},$$
where $\theta,a_1,a_2$ are all real parameters, subject to $\theta>0$
and $a_1a_2>0$. The set $\Lambda'$ is arguably more massive than
the lattice-cross $\Lambda$ of Corollary \ref{cor-1}. Nevertheless, if 
$\Gamma_\varepsilon$ is as in Corollary \ref{cor-1}, with $\varepsilon$ 
positive, it can be shown that $(\Gamma_\varepsilon,\Lambda')$ fails to be a 
Heisenberg uniqueness pair, no matter what positive values $\varepsilon$ and 
$\theta$ assume. Analogously, $(\Gamma_\varepsilon,\Lambda'')$ also fails
to be a Heisenberg uniqueness pair, for all $\varepsilon>0$ and $a_1a_2>0$
(but it can be shown that $(\Gamma_\varepsilon,\Lambda'\cup\Lambda'')$ is
a Heisenberg uniqueness pair, however).
This suggests that it is crucial that the points of the lattice-cross 
$\Lambda$ of Corollary \ref{cor-1} are located along the characteristic 
directions for the Klein-Gordon equation (the two axes).
\end{rem}

We need a result of algebraic nature.

\begin{lem}
Let $z_1,z_2\in\C$ be two points such that 
$$z_1-z_2=am\in a\Z,\quad \frac{1}{z_1}-\frac1{z_2}=b n\in b\Z,$$
for some positive reals $a,b$. Then, unless $z_1=z_2$, we have
$$z_1=\frac{am}{2}\bigg(1\pm\sqrt{1-\frac{4}{abmn}}\bigg),\quad 
z_2=z_1-am.$$ 
\label{lm-1}
\end{lem}

The proof is a simple exercise, and therefore omitted.

\begin{rem}
\label{rem-2}
Let us consider the singular measure $\mu=\delta_{u}-\delta_{v}$, where 
$$u=(u_1,1/u_1)\in\Gamma,\qquad v=(v_1,1/v_1)\in\Gamma.$$
Then
$$\widehat\mu(\xi)=\e^{\pi\imag(\xi_1u_1+\xi_2/u_1)}-
\e^{\pi\imag(\xi_1v_1+\xi_2/v_1)},$$
so that
$$\widehat\mu(\xi_1,0)=\e^{\pi\imag\xi_1u_1}-
\e^{\pi\imag\xi_1v_1},\qquad
\widehat\mu(0,\xi_2)=\e^{\pi\imag\xi_2/u_1}-
\e^{\pi\imag\xi_2/v_1}.$$
Suppose we try to achieve that
\begin{equation}
\widehat\mu(\alpha j,0)=\widehat\mu(0,\beta k)=0,\qquad j,k\in\Z,
\label{eq-zero}
\end{equation}
for some positive reals $\alpha,\beta$. 
We see that this amounts to  
$$\e^{\pi\imag\alpha u_1}=\e^{\pi\imag\alpha v_1}, \qquad 
\e^{\pi\imag\beta/u_1}=\e^{\pi\imag\beta/v_1},$$
which we rewrite in the form
$$u_1-v_1\in\frac{2}{\alpha}\Z,\qquad \frac{1}{u_1}-\frac{1}{v_1}\in
\frac{2}{\beta}\Z.$$
In view of Lemma \ref{lm-1}, there are plenty of such points $u_1,v_1\in\R$
with $u_1\neq v_1$, for any given $\alpha,\beta$. This shows that the 
 requirement that the measure $\mu$ be absolutely continuous with respect to 
arc length measure on $\Gamma$ is essential; without it, Theorem \ref{thm-1} 
would simply not be true.
\end{rem}

\section{Dynamics of a Gauss-type map}
\label{sec-dyn}

\noindent{\bf A Gauss-type map}. 
In order to prove our main theorem (Theorem \ref{thm-1}), we will need to 
study the invariant measures of a particular map.
We shall consider a map on the interval $]-1,1]$, which we think of as 
$\R/2\Z$ (topologically as well). 
The map in question is defined by $U(0)=0$ and
$$U(x)=\bigg\{-\frac{1}{x}\bigg\}_2,\qquad x\neq0,$$
where for real $t$, the expression $\{t\}_2\in]-1,1]$ is the unique number such
that $t-\{t\}_2\in2\Z$. The function $U$ is locally
strictly increasing and continuous, except for being interrupted by jumps.
The map $U:]-1,1]\to]-1,1]$ is associated with {\em continued fractions with 
even partial quotients} (see \cite{S1}, \cite{S2}, \cite{KL}, \cite{C}).
We see that, for $j=\pm1,\pm2,\pm3,\ldots$,
$$U(x)=-\frac{1}{x}+2j,\qquad \frac{1}{2j+1}<x\le\frac{1}{2j-1},$$
and hence $U$ maps the interval 
$]\frac{1}{2j+1},\frac{1}{2j-1}]$ onto $]-1,1]$ in a one-to-one fashion.
The derivative of $U$ is locally
$$U'(x)=\frac{1}{x^2}, \qquad x\in]-1,1]\setminus \frac{1}{2\Z+1}.$$
The point $1$ is a fixed point for $U$, and $U'(1^-)=U'(-1^+)=1$, which makes
$1$ is a {\em weakly repelling fixed point}. This means that when we iterate 
$U$, once we are close to $1$ (which is the same point as $-1$ in $\R/2\Z$), 
the successive iterates will remain near $1$ for a long time. If $x\in]-1,1]$
is rational, then after a finite number of steps, the $U$-iterate of $x$
is either $0$ or $1$ (see, for instance \cite{KL}). 
This illuminates why irrational numbers tend to spend a large portion of 
their $U$-orbits near $1$.
 
\medskip

\noindent\bf Invariant measures. \rm
If $\varphi$ is a continuous function on $\R/2\Z$ and $\nu$ is a bounded 
complex Borel measure on $]-1,1]$, then the integral
\begin{equation}
\int_{]-1,1]}\varphi(x)\,
\diff\nu(x)
\label{eq-100}
\end{equation}
is well-defined. However, the integral \eqref{eq-100} makes sense under weaker
assumptions on $\varphi$. Suppose $E$ is an open subset of $]-1,1]$ such that
the complement $]-1,1]\setminus E$ is countable, and that $\varphi$ is bounded
on $]-1,1]$ and continuous on $E$. Then \eqref{eq-100} makes sense for 
$\varphi$, and we call the function $\varphi$ {\em pseudo-continuous}.
We recall the familiar notion that a bounded complex Borel measure $\nu$ on 
$]-1,1]$ is $U$-{\em invariant} provided that
\begin{equation}
\int_{]-1,1]}\varphi(U(x))\,\diff\nu(x)=\int_{]-1,1]}\varphi(x)\,
\diff\nu(x)
\label{eq-101}
\end{equation}
holds for all pseudo-continuous test functions $\varphi$; it is easy to see
that $\varphi\circ U$ is pseudo-continuous if $\varphi$ is pseudo-continuous,
so that \eqref{eq-101} makes sense.
We shall reformulate this criterion in more concrete terms. First, we note 
that
$$\int_{]-1,1]\setminus\{0\}}\varphi(U(x))\,\diff\nu(x)=\sum_{j\in\Z^*}
\int_{]\frac{1}{2j+1},\frac{1}{2j-1}]}\varphi(U(x))\,\diff\nu(x)
=\sum_{j\in\Z^*}
\int_{]\frac{1}{2j+1},\frac{1}{2j-1}]}\varphi\bigg(-\frac1x+2j\bigg)\,
\diff\nu(x),$$
where $\Z^*=\Z\setminus\{0\}$, and that
$$\int_{]\frac{1}{2j+1},\frac{1}{2j-1}]}\varphi\bigg(-\frac1x+2j\bigg)\,
\diff\nu(x)=\int_{]-1,1]}\varphi(t)\,\diff\nu_j(t),$$
where
\begin{equation}
\diff\nu_j(t)=\diff\nu\bigg(
\frac{1}{2j-t}\bigg),\qquad -1<t\le1,
\label{eq-101.5}
\end{equation}
so that we have
$$\int_{]-1,1]\setminus\{0\}}\varphi(U(x))\,\diff\nu(x)
=\sum_{j\in\Z^*}\int_{]-1,1]}
\varphi(t)\,\diff\nu_j(t).$$
It follows that $\nu$ is $U$-invariant if and only if
\begin{equation}
\nu=\nu(\{0\})\delta_0+\sum_{j\in\Z^*}\nu_j.
\label{eq-102}
\end{equation}

More generally, given $\lambda\in\C$, we want to talk about 
$(U,\lambda)$-invariant measures, defined by the requirement that 
\begin{equation}
\int_{]-1,1]}\varphi(U(x))\,\diff\nu(x)=\lambda\int_{]-1,1]}\varphi(x)\,
\diff\nu(x)
\label{eq-1010}
\end{equation}
hold for all test functions $\varphi$; specifically, this means that
\begin{equation}
\lambda\nu=\nu(\{0\})\delta_0+\sum_{j\in\Z^*}\nu_j.
\label{eq-1020}
\end{equation}
It is easy to see that for $|\lambda|>1$, there are no $(U,\lambda)$-invariant
measures except for the zero measure.

\begin{prop}
Suppose $\nu$ is a bounded $(U,\lambda)$-invariant measure on $]-1,1]$, 
and write $\nu=\nu_a+\nu_s$, where $\nu_a$ is absolutely continuous, while 
$\nu_s$ is singular. 
Then $\nu_a$ and $\nu_s$ are also $(U,\lambda)$-invariant. Moreover, 
if $|\lambda|=1$, then  $|\nu|$, $|\nu_a|$, and $|\nu_s|$ are all $U$-invariant
measures.  
\label{prop-01}
\end{prop}

\begin{proof}
The relation \eqref{eq-1020} splits:
\begin{equation}
\nu_a=\lambda\sum_{j\in\Z^*}(\nu_j)_a,\qquad 
\nu_s=\lambda\bigg(\nu(\{0\})\delta_0+\sum_{j\in\Z^*}(\nu_j)_s\bigg),
\label{eq-103}
\end{equation}
where the subscripts $a$ and $s$ indicate the absolutely continuous and 
singular parts, respectively, of the measure in question.
We easily realize that $(\nu_j)_a=(\nu_a)_j$ and $(\nu_j)_s=(\nu_s)_j$,
so that \eqref{eq-103} expresses that $\nu_a$ and $\nu_s$ are both 
$U$-invariant.

Next, we suppose $|\lambda|=1$, and turn to the assertion that $|\nu|$ is 
$U$-invariant. Taking absolute values, we have
\begin{equation}
|\diff\nu(t)|\le|\nu(\{0\})|\diff\delta_0(t)+\sum_{j\in\Z^*}|\diff\nu_j(t)|,
\label{eq-2002}
\end{equation}
and so
\begin{multline*}
\int_{]-1,1]}|\diff\nu(t)|\le\sum_{j\in\Z^*}\int_{]-1,1]}|\diff\nu_j(t)|
=|\nu(\{0\})|+
\sum_{j\in\Z^*}\int_{]\frac{1}{2j+1},\frac{1}{2j-1}]}|\diff\nu(t)|\\
=|\nu(\{0\})|+\int_{]-1,1]\setminus\{0\}}|\diff\nu(t)|=
\int_{]-1,1]}|\diff\nu(t)|.
\end{multline*}
This is only possible if we have in fact equality in \eqref{eq-2002}:
$$|\diff\nu(t)|=|\nu(\{0\})|\diff\delta_0(t)
+\sum_{j\in\Z^*}|\diff\nu_j(t)|.$$
This relation expresses that $|\nu|$ is $U$-invariant; that $|\nu_a|$ and 
$|\nu_s|$ are $U$-invariant is a simple consequence of this fact.
\end{proof}
\medskip

\noindent \bf An unbounded smooth invariant measure. \rm
We now consider the positive unbounded smooth measure
$$\diff\omega(x)=\frac{\diff x}{1-x^2}.$$
The criterion \eqref{eq-1020} makes sense although $\omega$ is unbounded.
The following assertion was essentially found by Schweiger \cite{S1}.

\begin{prop}
The measure $\omega$ is $U$-invariant.
\end{prop}

We supply the simple proof.

\begin{proof}
We check that
$$\diff\omega_j(t)=\diff\omega\bigg(
\frac{1}{2j-t}\bigg)=\frac{\diff t}{(2j-t)^2-1},$$
and since
$$\sum_{j\in\Z^*}\frac1{(2j-t)^2-1}=
\frac{1}{2}\sum_{j\in\Z^*}\bigg(\frac1{2j-t-1}-\frac1{2j-t+1}\bigg)=
\frac12\bigg(\frac{1}{1+t}+\frac{1}{1-t}\bigg)=\frac{1}{1-t^2},$$
we find from \eqref{eq-1020} that $\omega$ is $U$-invariant.
\end{proof}

Schweiger \cite{S1} actually focused on the related map 
$|U|:[0,1]\to[0,1]$ given by $|U|(x)=|U(x)|$. He obtained the following 
basic result.

\begin{prop}
The measure $\omega$ is invariant also with respect to $|U|$. Moreover,
$|U|$ is ergodic, that is, if $E\subset[0,1]$ is a $|U|$-invariant set, then
either $\omega(E)=0$ or $\omega([0,1]\setminus E)=0$. 
\end{prop}
\medskip

\noindent\bf Consequences of Ergodic Theory. \rm 
The Birkhoff Ergodic Theorem -- in this setting of an unbounded invariant
ergodic measure \cite{Aa} -- states that if $\varphi$ is Borel measurable and
even with
$$\int_{-1}^{1}\frac{|\varphi(t)|}{1-t^2}\,\diff t<+\infty,$$
then 
$$\frac{1}{N}\sum_{k=0}^{N-1}\varphi(U^{\langle k\rangle}(t))\to0\quad
\text{as}\,\,\,N\to+\infty$$
almost everywhere on $]-1,1]$. Here, $U^{\langle k\rangle}$ stands for the 
$k$-th iterate of $U$. We observe that we do not need to know whether $U$ is
ergodic, just that $|U|$ is, if we use that $|U(-x)|=|U(x)|$. 
We pick $\varphi(t)=1-t^2$, and get:
\begin{equation}
\frac{1}{N}\sum_{k=0}^{N-1}\big(1-|U^{\langle k\rangle}(t)|^2\big)\to0\quad
\text{as}\,\,\,N\to+\infty
\label{eq-104}
\end{equation}
almost everywhere on $]-1,1]$. Suppose $\nu$ is a {\em positive, bounded, and
absolutely continuous} $U$-invariant measure on $]-1,1]$. 
By the $U$-invariance, we have
\begin{equation}
\int_{]-1,1]}\big(1-|U^{\langle k\rangle}(t)|^2\big)\,\diff\nu(t)
=\int_{]-1,1]}(1-t^2)\,\diff\nu(t),
\label{eq-104.5}
\end{equation}
and so
\begin{equation*}
\int_{]-1,1]}\frac{1}{N}\sum_{k=0}^{N-1}
\big(1-|U^{\langle k\rangle}(t)|^2\big)\,
\diff\nu(t)=\int_{]-1,1]}(1-t^2)\,\diff\nu(t).
\end{equation*}
By the Lebesgue dominated convergence theorem, it follows from
\eqref{eq-104} that 
\begin{equation*}
\int_{]-1,1]}\frac{1}{N}\sum_{k=0}^{N-1}
\big(1-|U^{\langle k\rangle}(t)|^2\big)\,
\diff\nu(t)\to0,\quad \text{as}\,\,\,N\to+\infty,
\end{equation*}
which combined with the \eqref{eq-104.5} leads to
$$\int_{]-1,1]}(1-t^2)\,\diff\nu(t)=0.$$
This is only possible if $\nu=0$.

We formalize this in a proposition.

\begin{prop}
Suppose $\lambda\in\C$ has $|\lambda|=1$, and that $\nu$ is an absolutely 
continuous bounded complex $(U,\lambda)$-invariant Borel measure on $]-1,1]$. 
Then $\nu=0$. 
\label{prop-2.4}
\end{prop}

\begin{proof}
By Proposition \ref{prop-01}, $|\nu|$ is a $U$-invariant measure. By the above
argument, $|\nu|=0$, and so $\nu=0$.
\end{proof}

\section{Extension of the trigonometric system}

\noindent\bf The trigonometric system. \rm
The trigonometric system $\{e_n(x)\}_{n\in{\mathbb Z}}$, with 
$e_n(x)=\e^{\pi\imag nx}$, is very successful in describing $2$-periodic 
functions on the line. Harald Bohr -- the brother of Niels Bohr, the physicist
-- developed over a number of years in the 1920s and 1930s the theory of almost
periodic functions based on more general real frequencies rather than the 
integer frequencies of the trigonometric system. 
\medskip

\noindent\bf An extension of the trigonometric system. \rm Here, 
we consider another 
extension of the trigonometric system, connected with the theory of 
composition operators. Let $\beta$ be a positive real parameter.
We introduce, for integers $n$, 
$$e^{\langle\beta\rangle}_n(x)=e_n\bigg(\frac{\beta}{x}\bigg)=
\e^{\pi\imag\beta n/x},$$
and note that these functions are bounded on the real line. 

After a dilation of the line, Theorem \ref{thm-1} is equivalent to the 
following statement.

\begin{thm}
As $n$ ranges over the integers, the functions $e_n(x)$ and 
$e^{\langle\beta\rangle}_n(x)$ form a weak-star-spanning system in 
$L^\infty(\R)$ if and only if $0<\beta\le1$.
\label{th-2}
\end{thm}

If $\mu$ is a positive bounded absolutely continuous Borel measure on $\R$,
then a bounded function in $L^\infty(\R)$ is automatically in
$L^p(\R,\mu)$ for $1<p<+\infty$, and the weak-star closure of a subspace
in $L^\infty(\R)$ is contained in the norm closure in $L^p(\R,\mu)$. 
We then have the following consequence of Theorem \ref{th-2}. The necessity
part just requires mimicking the corresponding argument involving harmonic
extensions in Section \ref{sec-4} below.

\begin{cor}
Suppose $1<p<+\infty$, and that $\diff\mu(x)=M(x)\diff x$, where $M(x)\ge0$
is Borel measurable, with
$$0<\int_{-\infty}^{+\infty}M(x)\,\diff x<+\infty.$$ 
Then, as $n$ ranges over the integers, the functions $e_n(x)$ and 
$e^{\langle\beta\rangle}_n(x)$ form a spanning 
system in $L^p(\R,\mu)$ provided that $0<\beta\le1$. If, in addition, 
$$\int_{-\infty}^{+\infty}\frac{\diff x}{(1+x^2)^{p/(p-1)}M(x)^{1/(p-1)}}
<+\infty,$$
the condition $0<\beta\le1$ is also necessary in order to have a spanning 
system.
\label{cor-101}
\end{cor}

\section{Necessity of the condition $0<\beta\le1$}
\label{sec-4}

\noindent\bf Harmonic extension. \rm
We extend the functions $e_n$ harmonically and boundedly to the upper half 
plane $\C_+=\{z\in\C:\,\,\im z>0\}$:
$$e_n(z)=\e^{\pi\imag n z},\qquad \im z\ge0, \,\,\,n\ge0,$$
while
$$e_n(z)=\e^{\pi\imag n \bar z},\qquad \im z\ge0,\,\,\,n<0.$$
Likewise, the harmonic extension of $e^{\langle\beta\rangle}_n$ is
$$e^{\langle\beta\rangle}_n(z)=\e^{\pi\imag\beta n/\bar z},\qquad \im z\ge0, 
\,\,\,n\ge0,$$
and
$$e^{\langle\beta\rangle}_n(z)=\e^{\pi\imag \beta n/z},\qquad 
\im z\ge0,\,\,\,n<0.$$
\smallskip

\noindent \bf Point separation. \rm
A general $L^\infty(\R)$ function is extended harmonically and boundedly to
$\C_+$ via the Poisson kernel; for each $z_0=x_0+\imag y_0\in\C_+$, the point 
evaluation functional $f\mapsto f(z_0)$ is given by
$$f(z_0)=\frac{1}{\pi}\int_{-\infty}^{+\infty}P(t,z_0)\,f(t)\,\diff t,\qquad
P(t,z_0)=\frac{y_0}{(x_0-t)^2+y_0^2},$$
where $t\mapsto P(t,z_0)$ is in $L^1(\R)$, and the functional is therefore
weak-star continuous on $L^\infty(\R)$. As we harmonically extend all the
functions in $L^\infty(\R)$, we get the space of all bounded harmonic functions
in $\C_+$. The bounded harmonic functions in $\C_+$ separate the points of
$\C_+$, so if we can find two points $z_1,z_2\in\C_+$ with $z_1\neq z_2$, 
such that 
\begin{equation}
e_n(z_1)=e_n(z_2),\qquad 
e^{\langle\beta\rangle}_n(z_1)=e^{\langle\beta\rangle}_n(z_2),
\label{eq-105}
\end{equation}
for all $n\in\Z$, then the linear span of $e_n,e^{\langle\beta\rangle}_n$,
cannot be weak-star dense in $L^\infty(\R)$. The condition \eqref{eq-105}
boils down to 
$$z_1-z_2\in2\Z,\qquad \frac{1}{z_1}-\frac{1}{z_2}\in\frac{2}{\beta}\Z,$$
where we may apply Lemma \ref{lm-1}, with $m=n=1$, $a=2$, and 
$b=2/\beta$. Assuming that $1<\beta<+\infty$,  we get that
$$z_1=1+\imag\sqrt{\beta^2-1},\qquad z_2=-1+\imag\sqrt{\beta^2-1},$$
are points in $\C_+$ satisfying \eqref{eq-105}. It follows that the requirement
$0<\beta\le1$ is necessary in Theorem \ref{th-2}.

\section{Periodic and inverted-periodic functions}

\noindent \bf Periodic and inverted-periodic functions. \rm
The weak-star closure in $L^\infty(\R)$ of the linear span of the functions 
$e_n(x)=\e^{\pi\imag nx}$, $n\in\Z$, equals $L^\infty_{2}(\R)$, the subspace
of $2$-periodic functions. Similarly, the weak-star closure of the 
linear span of the functions $e^{\langle\beta\rangle}(x)=e_n(\beta/x)$, 
$n\in\Z$, equals the subspace $L^\infty_{\langle\beta\rangle}(\R)$ of all 
functions $f\in L^\infty(\R)$ with $x\mapsto f(\beta/x)$ being $2$-periodic. 
Let us tacitly extend all functions in $L^\infty(\R)$ harmonically to $\C_+$ 
using the Poisson kernel. 
\medskip

\noindent\bf The intersection space. \rm 
Let us, for a moment, consider the intersection 
$$L^\infty_{2}(\R)\cap L^\infty_{\langle\beta\rangle}(\R).$$
We introduce $\group$ as the group of M\"obius transformations
preserving $\C_+$ generated by the translation $z\mapsto z+2$ and the
mapping $z\mapsto \beta z/(\beta-2z)$; then the elements of
$L^\infty_{2}(\R)\cap L^\infty_{\langle\beta\rangle}(\R)$ are precisely the 
functions in $L^\infty(\R)$ that are invariant under
$f\mapsto f\circ\gamma$, for $\gamma\in\group$. This situation is investigated 
in \S 11.4 of Beardon's book \cite{B}. For $0<\beta\le1$, the group
$\group$ is discrete and free (see, e. g., Gilman and Maskit \cite{GM}), and 
the fundamental domain (hyperbolic polygon) associated with 
$\C_+/\group$ is given by
\begin{equation}
{\mathfrak D}(\beta)=\bigg\{z\in\C_+:\,\,|\re z|<1,\,\,\,
\bigg|z-\frac{\beta}{2}\bigg|>\frac{\beta}{2},\,\,\,
\bigg|z+\frac{\beta}{2}\bigg|
>\frac{\beta}{2}\bigg\}.
\label{eq-FD}
\end{equation}
The domain ${\mathfrak D}(\beta)$ has a cusp at infinity and at the origin.
In addition, it has cusp(s) at $\pm1$ for $\beta=1$. 
For $0<\beta<1$, the fundamental domain has two boundary line segments 
$]-1,-\beta[$ and $]\beta,1[$, which is enough for $\C_+/\group$ to carry 
plenty of bounded harmonic (holomorphic as well) functions. 
A cusp is a removable singularity for a bounded harmonic function on 
$\C_+/\group$ (it is just an isolated removed point on the Riemann surface), 
which means that only constants are bounded and harmonic 
on $\C_+/\group$ for $\beta=1$. 
For $1<\beta<+\infty$, the group $\group$ is discrete if and only if 
\begin{equation}
\beta=\frac{1}{\cos^2(p\pi/(2q))}
\label{eq-beta}
\end{equation}
for some coprime positive integers $p,q$ with $p<q$ and $p\in\{1,2\}$. 
In case $p=1$, the fundamental domain is still given by \eqref{eq-FD}, while
for $p=2$ it is smaller, but retains two of the cusps. Anyway, under 
\eqref{eq-beta}, only cusps (two or three) occur in $\C_+/\group$, and all 
bounded harmonic functions are constant. In the remaining case, when 
\eqref{eq-beta} fails, the group $\group$ is non-discrete, and then 
every harmonic function which is invariant under $\group$ is necessarily 
constant. 

We gather some of the above observations in a proposition.

\begin{prop}
We have
$$L^\infty_{2}(\R)\cap L^\infty_{\langle\beta\rangle}(\R)=
\{\text{\rm constants}\}$$
if and only if $1\le\beta<+\infty$. Moreover, for $0<\beta<1$, 
$L^\infty_{2}(\R)\cap L^\infty_{\langle\beta\rangle}(\R)$ is 
infinite-dimensional.
\end{prop}
\medskip

\noindent\bf The sum space. \rm
Next, we turn to the study of the sum space. In order to obtain Theorem 
\ref{th-2}, we are to show that 
$$L^\infty_{2}(\R)+ L^\infty_{\langle\beta\rangle}(\R)$$
is weak-star dense in $L^\infty(\R)$ if and only if $0<\beta\le1$. 
In Section \ref{sec-4}, we saw that the sum fails to be be weak-star dense
for $1<\beta<+\infty$. In the sequel, we therefore {\em assume that} 
$0<\beta\le1$.
We now make a basic observation. Functions in 
$L^\infty_2(\R)$ may be prescribed freely on $]-1,1]$, but then they are
uniquely determined everywhere else, due to periodicity. Likewise, functions
in $L^\infty_{\langle\beta\rangle}(\R)$ are free on 
$\R\setminus]-\beta,\beta]$, and extend by ``periodicity'' everywhere else. 
This allows us to define operators $\Sop$, $\Tope_\beta$ as follows. 
The first operator,
$$\Sop:\,L^\infty(]-1,1])\to L^\infty(\R\setminus]-1,1])$$ 
is obtained by extending the function to be $2$-periodic on $\R$ and
then restricting the extended function to $\R\setminus]-1,1]$. 
The second operator,
$$\Tope_\beta:\,L^\infty(\R\setminus]-\beta,\beta])\to 
L^\infty(]-\beta,\beta])$$
is the analogous extension associated with the ``periodicity'' in
$L^\infty_{\langle\beta\rangle}(\R)$; in symbols, we may express it as
$$\Tope_\beta[f]=(\Sop[f\circ I_\beta])\circ I_\beta,$$
where $I_\beta(x)=-\beta/x$. 

Next, we agree on a useful convention.
For a Lebesgue measurable subset $X$ of the
real line of positive linear measure, we identify $L^\infty(X)$ with a 
weak-star closed subspace of $L^\infty(\R)$ by extending the functions to 
vanish on the complement $\R\setminus X$. 

\begin{lem}
If $\id$ is the identity operator. and if the operator 
$$\id-\Tope_\beta\Sop:\,L^\infty(]-1,1])\to 
L^\infty(]-1,1])$$
has weak-star dense range, then the sum space $L^\infty_{2}(\R)+
L^\infty_{\langle\beta\rangle}(\R)$ is weak-star dense in $L^\infty(\R)$. 
\label{prop-5.2}
\end{lem}

\begin{proof}
We write $\mathcal R$ for the range of the operator $\,\id-\Tope_\beta\Sop$.
Pick an arbitrary $F_2\in L^\infty(\R\setminus]-1,1])$, 
and ask of $F_1\in L^\infty(]-1,1])$ that 
$F_1=\Tope_\beta[F_2]+R$, where $R\in {\mathcal R}$. 
The set of all sums $F=F_1+F_2\in L^\infty(\R)$ we obtain in this fashion 
is denoted by $\mathcal F$. The following straightforward argument shows that
 $\mathcal F$ is weak-star dense in $L^\infty(\R)$. 
Suppose $K\in L^1(\R)$ has 
$$\langle F,K\rangle_\R=0$$
for all $F\in{\mathcal F}$.
We decompose $K=K_1+K_2\in L^1(\R)$, 
where 
$$K_1\in L^1(]-1,1])\quad\text{and}\quad K_2\in L^1(\R\setminus]-1,1]).$$ 
Then
$$0=\langle F,K\rangle_\R=\langle F_1,K_1\rangle_\R+\langle F_2,K_2\rangle_\R
=\langle\Tope_\beta[F_2],K_1\rangle_\R+\langle R,K_1\rangle_\R+
\langle F_2,K_2\rangle_\R.$$
We rewrite this in the form
$$\langle R,K_1\rangle_\R=-\langle F_2,K_2\rangle_\R-
\langle \Tope_\beta[F_2],K_1\rangle_\R.$$
Only the left hand side depends on $R\in{\mathcal R}$; by linearity, the 
only way this is possible is if $\langle R,K_1\rangle_\R=0$ for all 
$R\in{\mathcal R}$. But as ${\mathcal R}$ is dense we get that
$K_1=0$. The remaining relationship now reads
$$\langle F_2,K_2\rangle_\R=0.$$
As $F_2$ was arbitrary, we conclude that $K_2=0$ as well.

To finish the proof, we show that 
$${\mathcal F}\subset L^\infty_{2}(\R)+L^\infty_{\langle\beta\rangle}(\R).$$
For $F=F_1+F_2\in{\mathcal F}$ as above, the fact that $F_1-\Tope_\beta[F_2]
=R\in{\mathcal R}$ means that there exists a $g\in L^\infty(]-1,1])$ such that
\begin{equation}
g-\Tope_\beta\Sop[g]=(\id-\Tope_\beta\Sop)[g]=F_1-\Tope_\beta[F_2].
\label{eq-111}
\end{equation}
Also, let $h\in L^\infty(\R\setminus]-1,1])$ be given by
\begin{equation}
h=F_2-\Sop[g].
\label{eq-112}
\end{equation}
Then, by \eqref{eq-112}, 
\begin{equation}
F_2=h+\Sop[g],
\label{eq-112'}
\end{equation} 
so that
$$\Tope_\beta[F_2]=\Tope_\beta[h]+\Tope_\beta\Sop[g],$$
and if we combine this with \eqref{eq-111}, we get
\begin{equation}
F_1=\Tope_\beta[F_2]+g-\Tope_\beta\Sop[g]=\Tope_\beta[h]+\Tope_\beta\Sop[g]
+g-\Tope_\beta\Sop[g]=g+\Tope_\beta[h].
\label{eq-113}
\end{equation}
It follows from \eqref{eq-112'} and \eqref{eq-113} that
$$F=F_1+F_2=(g+\Tope_\beta[h])+(h+\Sop[g])=(g+\Sop[g])+(h+\Tope_\beta[h])
\in L^\infty_{2}(\R)+L^\infty_{\langle\beta\rangle}(\R).$$
The proof is complete.
\end{proof}

\medskip

\noindent\bf The operator $\Tope_\beta\Sop$. \rm
For $x\in\R$, let $\{x\}_2$ denote the number with $-1<\{x\}_2\le1$ and 
$x-\{x\}\in 2\Z$. Then since, for $\varphi\in L^\infty(]-1,1])$, 
$$\Sop[\varphi](x)=\varphi(\{x\}_2)\,\,1_{\R\setminus]-1,1]}(x),
\qquad x\in\R,$$
where $1_E$ denotes the characteristic function of the set $E\subset\R$,
we find that for $\psi\in L^\infty(\R\setminus]-\beta,\beta])$,
$$\Tope_\beta[\psi](x)=\psi\bigg(\frac{\beta}{\{\beta/x\}_2}\bigg)\,\,
1_{]-\beta,\beta]}(x),
\qquad x\in\R.$$
It follows that
\begin{equation}
\Tope_\beta\Sop[\varphi](x)=\varphi\bigg(\bigg\{
\frac{\beta}{\{\beta/x\}_2}\bigg\}_2\bigg)\,\,1_{E_\beta}(x),
\qquad x\in\R,
\label{eq-114}
\end{equation}
where
$$E_\beta=\bigg\{x\in]-\beta,\beta]\setminus\{0\}:
\,\,\frac{\beta}{\{\beta/x\}_2}\in\R\setminus]-1,1]\bigg\}.$$

\section{Analysis of a related composition operator}

\noindent\bf The Gauss-type map. \rm
Let $U_\beta$ be the mapping
$$U_\beta(x)=\{-\beta/x\}_2,\qquad x\neq0,$$
with $U_\beta(0)=0$. 
We consider the associated compressed composition operator 
$$\Cop_\beta:\,L^\infty(]-1,1])\to L^\infty(]-1,1])$$
given by
$$\Cop_\beta[f](x)=f(U_\beta(x))\,1_{]-\beta,\beta]}(x).$$
We quickly realize from \eqref{eq-114} that 
$$\Tope_\beta\Sop=\Cop_\beta^2,$$ 
and turn to analyzing $\Cop_\beta$. The identity
$$\id-\Tope_\beta\Sop=\id-\Cop_\beta^2=(\id+\Cop_\beta)(\id-\Cop_\beta)$$
shows that if $\id+\Cop_\beta$ and $\id-\Cop_\beta$ both have weak-star dense 
range, then $\id-\Cop_\beta^2$ has weak-star dense range as well. 
By elementary Functional Analysis, the operators $\id+\Cop_\beta$ and 
$\id-\Cop_\beta$ both have weak-star dense range if and only if for 
$\lambda=\pm1$, the (predual) adjoint
$$\lambda\id-\Cop_\beta^*:\,\,L^1(]-1,1])\to L^1(]-1,1])$$ 
has null kernel, that is, if the points $\pm1$ both fail to be eigenvalues of 
$\Cop_\beta^*$. The following result shows that this is the case for 
$0<\beta\le1$, making $\id-\Tope_\beta\Sop$ have weak-star dense range, and 
in view of Lemma \ref{prop-5.2}, then, 
$$L^\infty_2(\R)+L^\infty_{\langle\beta\rangle}(\R)$$
is weak-star dense in $L^\infty(\R)$. Given that we have verified the necessity
of the condition $0<\beta\le1$ in the context of Theorem \ref{th-2}, 
the rest of the assertion of Theorem \ref{th-2} follows.

\begin{prop}
For $0<\beta\le1$, the point spectrum $\sigma_p(\Cop^*_\beta)$ of 
$\Cop^*_\beta:\,L^1(]-1,1])\to L^1(]-1,1])$ is contained in the open unit 
disk $\D$. In particular, $\pm1$ are not eigenvalues of $\Cop^*_\beta$.
\end{prop}

\begin{proof}
It is clear that 
$$\sigma_p(\Cop^*_\beta)\subset\sigma(\Cop^*_\beta)\subset\bar\D,$$
where $\bar\D$ is the closed unit disk, so we just need to show that 
$\lambda\in\C$ with $|\lambda|=1$ cannot be eigenvalues.

The treatment of the cases $\beta=1$ and $0<\beta<1$ will be
different.
\medskip

\noindent\bf The case $\beta=1$. \rm
For $\beta=1$, $U_\beta=U_1=U$, the map we studied back in Section 
\ref{sec-dyn}, and $\Cop_1$ is
$$\Cop_1[f](x)=f(U(x))\,1_{]-1,1]}(x).$$
The defining relation for the (predual) adjoint 
$\Cop_1^*:L^1(]-1,1])\to L^1(]-1,1])$ is
$$\int_{]-1,1]}f(x)\,\Cop_1^*[g](x)\,\diff x=
\langle\Cop_1^*[g],f\rangle_\R=\langle g,\Cop_1[f]\rangle_\R=
\int_{]-1,1]}f(U(x))\,g(x)\,\diff x.$$
If $g$ is a nontrivial eigenfunction for $\Cop_1^*$ with eigenvalue $\lambda$,
then $\Cop_1^*[g]=\lambda g$, and so
$$\lambda\int_{]-1,1]}f(x)\,g(x)\,\diff x=
\int_{]-1,1]}f(U(x))\,g(x)\,\diff x;$$ 
this expresses that the absolutely continuous bounded measure 
$\diff\nu(x)=g(x)\diff x$ is a $(U,\lambda)$-invariant measure. 
By Proposition \ref{prop-2.4}, there are no bounded $(U,\lambda)$-invariant 
measures except the null measure, for $|\lambda|=1$. Consequently, 
$\lambda\in\C$ is not an eigenvalue of $\Cop_1^*$ for $|\lambda|=1$. 
\medskip

\noindent\bf The case $0<\beta<1$. \rm The same analysis reduces the problem 
to studying the bounded absolutely continuous measures $\nu$ on $]-1,1]$ with
\begin{equation}
\lambda\int_{]-1,1]}f(t)\,\diff\nu(t)=
\int_{]-1,1]}\Cop_\beta[f](t)\,\diff\nu(t)=
\int_{]-\beta,\beta]}f(U_\beta(t))\,\diff\nu(t) 
\label{eq-120}
\end{equation}
for all $f\in L^\infty(]-1,1])$, where $\lambda\in\C$ is fixed with 
$|\lambda|=1$. In more concrete terms, this amounts to
$$\lambda\,\diff\nu(t)=\sum_{j\in\Z^*}\diff\nu_j(t),\qquad t\in]-1,1],$$
where
$$\diff\nu_j(t)=\diff\nu\bigg(\frac{\beta}{2j-t}\bigg),\qquad t\in]-1,1].$$
Taking absolute values, we have, for $|\lambda|=1$,
$$|\diff\nu(t)|\le\sum_{j\in\Z^*}|\diff\nu_j(t)|,\qquad t\in]-1,1].$$
Integrating over $]-1,1]$, we find that
$$\int_{]-1,1]}|\diff\nu(t)|\le
\sum_{j\in\Z^*}\int_{]-1,1]}|\diff\nu_j(t)|=\int_{[-\beta,\beta]}
|\diff\nu(t)|,\qquad t\in]-1,1],$$
which is only possible if we have the equality
$$|\diff\nu(t)|=\sum_{j\in\Z^*}|\diff\nu_j(t)|,\qquad t\in]-1,1],$$
as well as
$$\diff\nu(t)=0,\qquad t\in ]-1,1]\setminus[-\beta,\beta].$$
If we iterate the relation \eqref{eq-120}, we get
\begin{equation}
\lambda^n\int_{]-1,1]}f(t)\,\diff\nu(t)=
\int_{]-1,1]}\Cop_\beta^n[f](t)\,\diff\nu(t)=
\int_{E_\beta(n)}f(U^{\langle n\rangle}_\beta(t))\,\diff\nu(t), 
\label{eq-121}
\end{equation}
where set $E_\beta(n)$ is given by 
$$E_\beta(n)=\big\{t\in]-1,1]:\,U^{\langle k\rangle}_\beta(t)\in[-\beta,\beta]
\,\,\,\text{for}\,\,\,k=0,\ldots,n-1\big\}.$$
A repetition of the above argument involving $U^{\langle n\rangle}_\beta$ in 
place of $U_\beta$ shows that if $|\lambda|=1$, then
$$\diff\nu(t)=0,\qquad t\in ]-1,1]\setminus E_\beta(n).$$
As $n\to+\infty$, the set $E_\beta(n)$ shrinks down to 
$$E_\beta(\infty)=
\big\{t\in]-1,1]:\,U^{\langle k\rangle}_\beta(t)\in[-\beta,\beta]
\,\,\,\text{for}\,\,\,k=0,1,2,3,\ldots\big\}.$$
This final set $E_\beta(\infty)$ is $U_\beta$-invariant, and it is not hard 
to show that it must have zero length. But the measure $\nu$ vanishes 
everywhere else, and being absolutely continuous, it must be the zero measure.
In particular, $\lambda\in\C$ with $|\lambda|=1$ cannot be eigenvalues of 
$\Cop_\beta^*$.
The proof is complete.
\end{proof}
\medskip

\noindent \bf A remark on model subspaces. \rm Given an inner function 
$\Theta$ in the upper half plane $\C_+$, one considers the model subspaces 
$K_\Theta(\C_+)=H^2(\C_+)\ominus\Theta H^2(\C_+)$. Uniqueness sets for 
model subspaces have been studied recently by Makarov and Poltoratski
\cite{MP}, and the injectivity of the Toeplitz operator with symbol
$\bar\Theta B_\Lambda$ is equivalent to $\Lambda\subset\C_+$ being a 
uniqueness set. In our setting, we use mainly that the operators 
$\lambda\id-\Cop_\beta^*$ are injective for $\lambda=\pm1$, so apparently 
these operators are analogous to the Toeplitz operators from the model 
subspace case. 

\section{Applications and open problems}

\noindent\bf An application to BMO. \rm
Let ${\text BMOA}(\C_+)$ be the (weak-star closed) subspace of the space
$\text{BMO}(\R)$ 
consisting of functions whose Poisson extensions to $\C_+$ are holomorphic in
$\C_+$. We recall that $\text{BMO}(\R)$ denotes the space of functions with 
bounded mean oscillation. The Cauchy-Szeg\"o (analytic) projection 
$${\mathbf P}:\,L^\infty(\R)\to{\text{BMOA}}(\C_+)$$
is bounded and surjective. We observe that if $\mathcal X$ is a linear 
subspace of $L^\infty(\R)$ which is weak-star dense, then 
${\mathbf P}({\mathcal X})$ is dense in ${\text{BMOA}}(\C_+)$. Moreover, 
let ${\text{BMOA}}_2(\C_+)$ be the subspace of ${\text BMOA}(\C_+)$ of 
functions invariant under $z\mapsto z+2$, and let 
${\text{BMOA}}_{\langle\beta\rangle}(\C_+)$ be the subspace of 
${\text{BMOA}}(\C_+)$ of functions invariant under 
$z\mapsto\beta z/(\beta-2z)$. We quickly check that ${\mathbf P}$ maps 
$L^\infty_2(\R)\to{\text{BMOA}}_2(\C_+)$ and $L^\infty_{\langle\beta\rangle}
(\R)\to{\text{BMOA}}_{\langle\beta\rangle}(\C_+)$.

In view of our main theorem (Theorem \ref{th-2}), we have the following.

\begin{cor} The sum ${\text{\rm BMOA}}_2(\C_+)+
{\text{\rm BMOA}}_{\langle\beta\rangle}(\C_+)$
is weak-star dense in $\text{\rm BMOA}(\C_+)$ if and only if $0<\beta\le1$. 
In other words, the functions
$$e_{n}(x)=\e^{\pi\imag nx},\quad e^{\langle\beta\rangle}_{-n}(x)=
\e^{-\pi\beta\imag n/x},\qquad n=0,1,2,3,\ldots,$$
span a weak-star dense subspace of $\text{\rm BMOA}(\C_+)$ if and only if 
$0<\beta\le1$. 
\label{cor-102}
\end{cor}

By the M\"obius invariance of BMO, we may transfer this 
result to the setting of the unit disk, and answer Problem 2 of Matheson and 
Stessin \cite{MS} in the affirmative. 
\medskip

\noindent\bf Four open problems. \rm 
(a) Suppose in the context of Theorem \ref{thm-1} we consider a lattice-cross
$$\Lambda=((\alpha\Z+\{\theta\})\times\{0\})\cup(\{0\}\times\beta\Z),$$
where $\theta\in\R$ is fixed. It seems that Theorem \ref{thm-1} should remain
true with this new $\Lambda$, with only moderate modifications in the proof.
But what happens if the lattice-cross is less regular, that is, if the two 
spacings $\alpha$ and $\beta$ are allowed to fluctuate a bit along the cross? 

\noindent{(b)} In the context of Corollary \ref{cor-101},  as $n$ ranges over
the integers, do the functions $e_n(x)$ and $e^{\langle\beta\rangle}_n(x)$
form a spanning system in $L^p(\R,\mu)$ for all $0<\beta<+\infty$ provided
that
$$\int_{-\infty}^{+\infty}\frac{\diff x}{(1+x^2)^{p/(p-1)}M(x)^{1/(p-1)}}
=+\infty?$$

\noindent{(c)} Let $H^\infty_2(\C_+)$ denote the 
(weak-star closed) subspace of $L^\infty_2(\R)$ consisting of those functions
whose Poisson extensions to the upper half plane $\C_+$ are analytic. 
Analogously, let $H^\infty_{\langle\beta\rangle}(\C_+)$ be the 
(weak-star closed) subspace of $L^\infty_{\langle\beta\rangle}(\R)$ 
consisting of those functions whose Poisson extensions to the upper half 
plane $\C_+$ are analytic. Is the sum 
$$H^\infty_2(\C_+)+H^\infty_{\langle\beta\rangle}(\C_+)$$
weak-star dense in $H^\infty(\C_+)$ for $0<\beta\le1$? This does not seem to 
follow from our Theorem \ref{th-2}, and, if answered in the affirmative (for
$0<\beta<1$), would solve Problem 1 of \cite{MS}. 

\noindent{(d)} In the context of Corollary \ref{cor-1}, let the bounded 
Borel measure $\mu$ be supported
on the hyperbola $x_1x_2=\varepsilon$, and let $\Lambda$ be the lattice-cross
given by the positive parameters $\alpha,\beta$. If $\mu$ is absolutely 
continuous with respect to arc length measure on the hyperbola, an argument
involving curvature considerations shows that
$$\widehat\mu(\xi)\to0\quad\text{as}\,\,\,|\xi|\to+\infty,$$
and this in a sense expresses the absence of point masses in $\mu$. Moreover,
$\widehat\mu$ solves the Klein-Gordon equation
$\partial_1\partial_2\widehat\mu+\varepsilon\pi^2\widehat\mu=0$. It would
be desirable to remove to the extent possible the Fourier analysis ingredient
in Corollary \ref{cor-1}. Let $u$ be a bounded continuous complex-valued
function on $\R^2$ with $u(\xi)\to0$ as $|\xi|\to+\infty$. Suppose, in 
addition, that $\partial_1\partial_2u+\varepsilon\pi^2u=0$ holds in the
sense of distribution theory. The problem: is the lattice-cross $\Lambda$, with
 $\alpha\beta\le1/|\varepsilon|$, a uniqueness set for $u$ (that is, 
$u|_\Lambda=0\implies u=0$)?

\medskip

\noindent\bf Acknowledgements. \rm We thank Dani Blasi Babot, Francesco 
Cellarosi, Joaquim Ortega-Cerd\`a, Eero Saksman, Masha Saprykina, Yakov 
Sinai, and Serguei Shimorin for helpful comments and references.


\end{document}